\documentclass[11pt]{amsart}
\usepackage{amssymb,amsmath,hyperref,tikz,graphicx,float}
\usepackage{soul}
\usepackage{breqn} 

\newtheorem{theorem}{Theorem}
\newtheorem*{theorem*}{Theorem}

\newtheorem{conjecture}[theorem]{Conjecture}
\newtheorem*{conjecture*}{Conjecture}
\newtheorem{corollary}[theorem]{Corollary}
\newtheorem*{corollary*}{Corollary}
\newtheorem{definition}[theorem]{Definition}
\newtheorem{example}[theorem]{Example}

\newtheorem{exercise}[theorem]{Exercise}
\newtheorem{fact}[theorem]{Fact}
\newtheorem{lemma}[theorem]{Lemma}
\newtheorem{problem}[theorem]{Problem}
\newtheorem{proposition}[theorem]{Proposition}
\newtheorem*{proposition*}{Proposition}
\newtheorem{question}[theorem]{Question}
\newtheorem{remark}[theorem]{Remark}

\newcommand{\bcon}{\begin{conjecture}}
\newcommand{\econ}{\end{conjecture}}
\newcommand{\bcor}{\begin{corollary}}
\newcommand{\ecor}{\end{corollary}}
\newcommand{\bdf}{\begin{definition}}
\newcommand{\edf}{\end{definition}}
\newcommand{\benu}{\begin{enumerate}}
\newcommand{\eenu}{\end{enumerate}}
\newcommand{\beq}{\begin{equation}}
\newcommand{\eeq}{\end{equation}}
\newcommand{\bexa}{\begin{example}}
\newcommand{\eexa}{\end{example}}
\newcommand{\bexe}{\begin{exercise}}
\newcommand{\eexe}{\end{exercise}}
\newcommand{\bfac}{\begin{fact}}
\newcommand{\efac}{\end{fact}}
\newcommand{\bite}{\begin{itemize}}
\newcommand{\eite}{\end{itemize}}
\newcommand{\blem}{\begin{lemma}}
\newcommand{\elem}{\end{lemma}}
\newcommand{\bmat}{\begin{matrix}}
\newcommand{\emat}{\end{matrix}}
\newcommand{\bprb}{\begin{problem}}
\newcommand{\eprb}{\end{problem}}
\newcommand{\bpro}{\begin{proposition}}
\newcommand{\epro}{\end{proposition}}

\newcommand{\bque}{\begin{question}}
\newcommand{\eque}{\end{question}}
\newcommand{\brem}{\begin{remark}}
\newcommand{\erem}{\end{remark}}
\newcommand{\bthm}{\begin{theorem}}
\newcommand{\ethm}{\end{theorem}}

\newtheorem*{namedtheorem}{\theoremname}
\newcommand{\theoremname}{testing}
\newenvironment{named}[1]{\renewcommand{\theoremname}{#1}\begin{namedtheorem}}{\end{namedtheorem}}

\newcommand{\bpr}{\begin{proof}}
\newcommand{\epr}{\end{proof}}

\newcommand{\lb}{\label}
\newcommand{\la}{\langle}
\newcommand{\ra}{\rangle}
\newcommand{\comment}[1]{\,}
\newcommand{\wh}{\widehat}
\newcommand{\cal}{\mathcal}
\newcommand{\p}{\partial}

\newcommand{\Z}{\mathbb Z}
\newcommand{\Q}{\mathbb Q}
\newcommand{\R}{\mathbb R}

\newcommand{\ve}{\varepsilon}

\newcommand{\vs}{\vspace*{.1in}}
\newcommand{\pmo}{{\pm 1}}

\newcommand{\diag}[2]{\parbox{#2}{\includegraphics[width=#2]{#1}}} 

\newcommand{\KP}[1]{
  \begin{tikzpicture}[baseline=-\dimexpr\fontdimen22\textfont2\relax]
  #1
  \end{tikzpicture}
}
\newcommand{\KPC}{
  \KP{\filldraw[color=black, fill=none, thick] circle (0.18);}
}
\newcommand{\KPX}{
  \KP{
    \draw[color=black, thick] (-0.3,-0.3) -- (0.3,0.3);
    \draw[color=black, thick] (-0.3,0.3) -- (-0.05,0.05);
    \draw[color=black, thick] (0.05,-0.05) -- (0.3,-0.3);
  }
}
\newcommand{\KPB}{
  \KP{
    \draw[color=black, thick] (-0.3,0.3) .. controls (0,-0.02) .. (0.3,0.3);
    \draw[color=black, thick] (-0.3,-0.3) .. controls (0,0.02) .. (0.3,-0.3);
  }
}
\newcommand{\KPA}{
  \KP{
    \draw[color=black, thick] (-0.3,-0.3) .. controls (0.02,0) .. (-0.3,0.3);
    \draw[color=black, thick] (0.3,-0.3) .. controls (-0.02,0) .. (0.3,0.3);
  }
}

\newcommand{\KPAKink}{
  \KP{
    \draw[color=black, thick] (-0.3,-0.3) .. controls (-0.15,0) .. (-0.3,0.3);
    \draw[color=black, thick] (0.3,-0.3) -- (0.16,-0.08);
    \draw[color=black, thick] (0.3,0.3) .. controls (-0.2,-0.8) and (-0.2,0.5) .. (0.1,0.05);
    {}
  }
}

\newcommand{\TwithCross}{
  \KP{
    \filldraw[color=black, fill=none, thick] circle (0.25);
    \put(-0.05,-0.05){$T$}
    \draw[color=black, thick] (-0.3,-0.3) -- (-0.19,-0.19);
    \draw[color=black, thick] (-0.3,0.3) -- (-0.19,0.19);
    \draw[color=black, thick] (0.6,-0.3) -- (0.52,-0.1);
    \draw[color=black, thick] (0.19,0.19) .. controls (0.35,0.3) and (0.4,0.3) .. (0.48,0.1);
    \draw[color=black, thick] (0.6,0.3) .. controls (0.45,-0.3) and (0.35,-0.3) .. (0.19,-0.19);
    {}
    {}
  }
}

\title{Tangle Equations, the Jones conjecture, slopes of surfaces in tangle complements, and $q$-deformed rationals
}

\author{Adam S. Sikora}

\thanks{The author was partially supported by the Simons Foundation grant 957582.}


\begin{document}

\thispagestyle{empty}

\begin{abstract} 
We study systems of $2$-tangle equations
$$\begin{cases}
 N(X+T_1)=L_1\\
 N(X+T_2)=L_2.
\end{cases}$$
which play an important role in the analysis of enzyme actions on DNA strands.

We show that every system of framed tangle equations has at most one framed rational solution. 
Furthermore, we show that the Jones Unknot conjecture implies that if a system of tangle equations has a rational solution then that solution is unique among all $2$-tangles. This result potentially opens a door to a purely topological disproof of the Jones Unknot conjecture.

We introduce the notion of the Kauffman bracket ratio $\{T\}_q\in \Q(q)$ of any $2$-tangle $T$ and we 
conjecture that for $q=1$ it is the slope of meridionally incompressible surfaces in $D^3-T$. We prove that conjecture for algebraic $T$. We also prove that for rational $T$, the brackets $\{T\}_q$ coincide with the $q$-rationals of \cite{MO}. 

Additionally,  we relate systems of tangle equations to the Cosmetic Surgery Conjecture and the Nugatory Crossing Conjecture. 

\end{abstract}
\address{244 Math Bldg, University at Buffalo, SUNY, Buffalo, NY 14260}
\email{asikora@buffalo.edu}

\pagestyle{myheadings}

\maketitle

\markboth{\hfil{\sc ADAM S. SIKORA}\hfil}
{\hfil{\sc Tangle Equations, the Jones conjecture, $q$-deformed rationals.}\hfil}

\newcommand{\TaddNumer}[2]{
  \KP{
    \filldraw[color=black, fill=none, thick] circle (0.28);
    \filldraw[color=black, fill=none, thick] (0.8,0) circle (0.28);
    \put(-0.07,-0.05){#1} 
    \put(0.235,-0.05){#2} 
    {}
    {}
    \draw[color=black, thick] (-0.19,-0.19) .. controls (-0.7, -0.7) and (1.5,-0.7) .. (0.99,-0.19);
    {}
    {}
    \draw[color=black, thick] (-0.19,0.19) .. controls (-0.7, 0.7) and (1.5,0.7) .. (0.99,0.19);
    \draw[color=black, thick] (0.2,0.19) .. controls (0.4, 0.3)  .. (0.60,0.19);
    \draw[color=black, thick] (0.2,-0.19) .. controls (0.4,-0.3)  .. (0.60,-0.19);
  }
}


\section{Introduction}


\subsection{Systems of tangle equations}
A \underline{system of tangle equations} has the form:
\beq\label{e-tangle0}
\TaddNumer{$X$}{$T_1$} =L_1\quad \text{and}\quad 
\TaddNumer{$X$}{$T_2$} =L_2,
\eeq
where $T_1\ne T_2$ are rational $2$-tangles and $L_1,L_2$ are links and $X$ is an unknown $2$-tangle. (From now on 2-tangles will be referred to as ``tangles''.) Such systems play an important role in the analysis of entanglement of DNA molecules, see Sec. \ref{s-systangleEq}. 
We study them in the settings of both unframed and framed links, by combining topological methods (surgery theory on double branched covers of links and tangles) and algebraic methods involving the Jones polynomial and the Kauffman bracket. 

We prove foundational results in Subsections \ref{s-tangleEq}-\ref{s-systangleEq}. In particular, we show that every tangle equation has infinitely many solutions when $L$ is not the unlink of two components. We also exhibit systems \eqref{e-tangle0} which have an arbitrarily large numbers of solutions.

The \underline{Kauffman bracket} of a framed 2-tangle $T$ is
$[T]=\begin{pmatrix} [T]_0\\ [T]_\infty\end{pmatrix},$ where 
$[T]_0,[T]_\infty\in \Z[A^{\pmo}]$ are determined by applying the Kauffman bracket skein relations reducing $T$ to a linear combination
$$T=[T]_0\cdot \KPB + [T]_\infty\cdot \KPA,$$
cf. \cite{EKT, LR, TS1}. 

A rational 
solution of  (\ref{e-tangle0}) is one given by a rational tangle.
We prove

\bpro[Proof in Sec. \ref{s.uniqueKB}]\label{p.rat-uniqe-i}
Every framed system (\ref{e-tangle0}) has at most one framed rational solution. 
\epro

We propose more generally:
\bcon\label{c.rationalunique}
If a framed system (\ref{e-tangle0}) has a rational solution $X$ then it is unique (among all solutions $X$).
\econ

The Jones Unknot conjecture, henceforth called the Jones Conjecture, states that the Jones polynomial distinguishes all non-trivial knots from the trivial one, \cite{Jo}.
One of the main results of this paper relates the above two conjectures:

\bthm[Proof in Sec. \ref{s-JC-for-tan}]\label{JC-UTESi}
The Jones conjecture implies Conjecture \ref{c.rationalunique}.
\ethm

A stronger version of this statement, formulated in Sec. \ref{s-JC-for-tan}, shows that because the Jones Conjecture holds up to 24 crossings, \cite{TS2},  Conjecture \ref{c.rationalunique} holds for tangles up to $11$ crossings. 

The topological meaning of the Jones polynomial is rather obscure, despite much research devoted to it in the last decades. Hence, we find this purely topological consequence of the Jones conjecture intriguing. In particular, this result opens a door to disproving the Jones conjecture through purely topological methods.

We will show in Proposition \ref{p.solnum} that there is no upper bound on the number of solutions of systems (\ref{e-tangle0}) (framed and unframed ones). 

Conjecture \ref{c.rationalunique} does not hold for unframed systems.
However, every unframed system has at most two rational solutions, \cite{ES}.
Conjecture \ref{c.rationalunique} suggests the following unframed version of it:

\bcon
If an unframed system (\ref{e-tangle0}) has a rational solution $X$ then it its every solution is rational.
\econ

The importance of this conjecture stems from the fact that a counterexample to it may lift to a framed counterexample of Conjecture \ref{c.rationalunique} thus disproving the Jones conjecture.

\subsection{Approach through the surgery theory}

Taking double branched covers of tangles and of links translates equations (\ref{e-tangle0}) into the language of surgery theory of $3$-manifolds. Apart from the Kauffman bracket, this is the second main source of methods utilized in this paper. In particular, we show

\def\pLL{If $L_1=L_2$ is either a nontrivial rational knot or a non-trivial torus knot then System (\ref{e-tangle0}) has no algebraic tangle (framed or unframed) solutions for any $T_1,T_2.$}

\bpro[Proof in Sec. \ref{s.surgery}.]\label{p.LL-i}
\pLL
\epro

(This result settles in particular Conjecture \ref{c.rationalunique} when $L_1=L_2$ is a nontrivial rational knot or a torus knot.) We will relate systems of tangle equations to Cosmetic Surgery Conjecture in Sec. \ref{s.surgery}.

\subsection{Connection with $q$-rationals and slopes of essential surfaces}

\bthm[Proof in Sec. \ref{s-q-rationals}.]\label{t.q-qi}
For any tangle $T$, 
$$\{T\}_q:=A^{-2} [T]_\infty/[T]_0$$
is well defined (i.e. not $0/0$) and independent of the framing of $T$.
Furthermore, 
$$\{T\}_q\in \wh{\Q(q)}=\Q(q)\cup \{\infty\},$$
where $q=-A^4.$
\ethm

We call $\{T\}_q$ the \underline{Kauffman bracket ratio} of $T$ or the \underline{KB ratio} for short.

\bthm[Proof in Sec. \ref{s-q-rationals}.]\label{t.qrati}
For any $x\in \wh{\Q}$, 
$\{\la x\ra\}_q$ is the $q$-rational of \cite{MO}.
\ethm

$q$-rationals are $q$-deformations of rational numbers, in the sense that $\{\la x\ra\}_q$ evaluates to $x$ at $q=1.$
(They were independently discovered by us in the first version of this paper on arXiv.)
Interestingly much of the theory of continued fractions extends to these $q$-rational numbers.
We believe that our approach to $q$-rational numbers provides a useful new geometric intuition. 

Finally, we relate the Kauffman bracket ratios to slopes of dividing surfaces in tangle complements in Sec. \ref{s.slope}.  A surface $S$ properly embedded in the complement of a tangle, $B^3-T$, is \underline{meridionally incompressible} if every $2$-disk $D^2$ in $B^3$ with boundary in $S$ intersecting $T$ transversely precisely once can be deformed to a disk in $S\cup T$. 
A surface $S$ is \underline{$m$-essential} (for ``meridionally essential'') if it is incompressible, meridionally incompressible and not boundary-parallel in $D^3 \smallsetminus T,$ see \cite{Oz1, Oz2}. 

We say that a loop $\gamma$ in $S^2_4$ is \underline{dividing} if it separates the punctures of $S^2_4$ into two groups of two. Such loops are classified by their slope $s(\gamma)\in \wh \Q=\Q\cup \{\infty\}$ which can be defined by lifting $\gamma$ to the torus being the double cover of $S^2$ branched over the four ends of $T$, see Sec. \ref{s.slope}. 

If the boundary of an $m$-essential surface contains a dividing loop $\gamma$ we call $S$ \underline{dividing} and we call $s(\gamma)$ the \underline{slope} $s(S)$ of $S$. (Then it is easy to see that all other dividing loops of $\p S$ must be parallel.) 
This definition is inspired by that in \cite{Oz2}, where the slopes are defined for algebraic tangles.

\bcon\label{c.slopeeq}
Every dividing $m$-essential surface in the complement of any tangle $T$ has slope $\{T\}_1$.
\econ

For the above reason, we call $\{T\}_1\in \wh\Q$  the \underline{algebraic slope} of $T$.
This conjecture is related to \cite[Question 3.4]{Oz2} which asks more broadly if the slopes of such surfaces are determined by $T$.

\bthm[Proof in Sec. \ref{s.slope}]\label{t.slope-i}
Our algebraic slope coincides with Ozawa's slope for algebraic tangles.
In other words, the above conjecture holds for algebraic tangles. 
\ethm

 {}

\textbf{Acknowledgements:} The author thanks V. Ovsienko and M. Ozawa for helpful discussions.

\newcommand{\Taddition}[2]{
  \KP{
    \filldraw[color=black, fill=none, thick] circle (0.28);
    \filldraw[color=black, fill=none, thick] (0.8,0) circle (0.28);
    \put(-0.07,-0.05){#1}
    \put(0.235,-0.05){#2}
    \draw[color=black, thick] (-0.3,-0.3) -- (-0.19,-0.19);
    \draw[color=black, thick] (-0.3,0.3) -- (-0.19,0.19);
    \draw[color=black, thick] (1.1,-0.3) -- (.99,-0.19);
    \draw[color=black, thick] (1.1,0.3) -- (0.99,0.19);
    \draw[color=black, thick] (0.2,0.19) .. controls (0.4, 0.3)  .. (0.60,0.19);
    \draw[color=black, thick] (0.2,-0.19) .. controls (0.4,-0.3)  .. (0.60,-0.19);
    {}
    {}
  }
}


\section{Preliminaries}

\subsection{Rational and Algebraic Tangles and Links}
Throughout, $2$-tangles will be called \underline{tangles} for brevity. They are proper tame embeddings of 
$1$-manifolds into $D^3=D^2\times (-1,1)$ with 
ends at NE, SE, SW, NW points of $\p D^2,$ where $D^2$ is a compact disk identified with $D^2\times \{0\}.$
Tangles are considered up to isotopies fixing $\p D^2$.
A \underline{$2$-string tangle} is one which consists of two strings only (no loops).
The $\la-1\ra , \la0\ra , \la1\ra $ and $\la\infty\ra $ tangles and the \underline{tangle addition} are depicted in Fig. \ref{f-basictangles}. 
\begin{figure}[h]
\KP{
   \filldraw[color=black, fill=none, thick] (0,0) circle (0.41);
    \draw[color=black, thick] (-0.3,0.3) -- (0.3,-0.3);
    \draw[color=black, thick] (-0.3,-0.3) -- (-0.05,-0.05);
    \draw[color=black, thick] (0.05,0.05) -- (0.3,0.3);
    {}
  }\quad
  \KP{
   \filldraw[color=black, fill=none, thick] (0,0) circle (0.41);
    \draw[color=black, thick] (-0.3,0.3) .. controls (0,-0.02) .. (0.3,0.3);
    \draw[color=black, thick] (-0.3,-0.3) .. controls (0,0.02) .. (0.3,-0.3);
  }\quad
   \KP{
   \filldraw[color=black, fill=none, thick] (0,0) circle (0.41);
    \draw[color=black, thick] (-0.3,-0.3) -- (0.3,0.3);
    \draw[color=black, thick] (-0.3,0.3) -- (-0.05,0.05);
    \draw[color=black, thick] (0.05,-0.05) -- (0.3,-0.3);
  }\quad
 \KP{
  \filldraw[color=black, fill=none, thick] circle (0.41);
  \draw[color=black, thick] (-0.3,-0.3) .. controls (0.02,0) .. (-0.3,0.3);
   \draw[color=black, thick] (0.3,-0.3) .. controls (-0.02,0) .. (0.3,0.3);
    {}
  }\quad
  \Taddition{$T$}{$T'$}
  {}
  {}
\caption{The $-1, 0, 1$ and $\infty$ tangles and the tangle addition $T+T'$. (We follow here Conway's notation, \cite{Co}. Kauffman's and his collaborators' papers use opposite signs, eg. \cite{KL}.)} 
\lb{f-basictangles}
\end{figure}

The result of adding $n$ tangles $\la 1\ra $ (respectively: $\la -1\ra $) together is denoted by $\la n\ra $ (respectively $\la -n\ra $), for $n=1,2,3...$
These tangles, together with $\la 0\ra$, are called \underline{integral.}

The \underline{mirror image}  $-T$ of $T$ 
is obtained by switching all crossings of $T$.
The \underline{tangle rotation} $R(T)$ 
is the $90^0$ clockwise rotation and the \underline{tangle inversion} 
is the tangle rotation followed by the mirror image. 

All tangles obtained from the integral ones by the operations of addition and rotation are called \underline{algebraic}. (This class is closed under the mirror image and inversion.) Among them are rational tangles defined as follows:

By $\la a_n,..,a_1\ra $ we denote the tangle obtained from $\la 0\ra $ by adding $\la a_1\ra $ followed by the inversion, then by adding $\la a_2\ra $ followed by the inversion, and so on, until this construction is finished by adding $a_n$ at the end, as in Figure \ref{f-rattan-closures}(left). Tangles of this form, for $a_1,...,a_n\in \Z$, are called \underline{rational} because as observed by Conway \cite{Co} (and proved in \cite{Mo2,BZ,KL})
$$\la a_n,..,a_1\ra  \to a_n+\frac{1}{a_{n-1}+\frac{1}{\ddots +\frac{1}{a_1}}}$$
defines a bijection between rational tangles (up to tangle isotopy) and 
\mbox{$\widehat \Q=\Q\cup \{\infty\}$.}
Under this bijection, tangle inversion corresponds to the fraction inversion, \mbox{$x\to 1/x,$} and the mirror image operation corresponds to the negation, $x\to -x.$

This association does not preserve addition.

By applying the numerator or the denumerator closure (defined in Fig. \ref{f-rattan-closures}) to a rational tangle we obtain a \underline{rational link}, also referred to as a $2$-bridge link or $4$-plat.

\begin{figure}[h]
\centerline{\diag{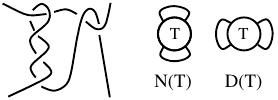}{2in}}
\caption{Rational Tangle $\la -2,-3,2\ra$, and the numerator, and the denumerator closures}
\lb{f-rattan-closures}
\end{figure}

\subsection{Tangle Equations}
\lb{s-tangleEq}

A \underline{tangle equation} has the form
\beq\lb{e-tangle}
 N(X+T)=L.
\eeq
where a rational tangle $T$ and a link $L$ are given and the tangle $X$ is unknown.
Finding all solutions of \eqref{e-tangle} is usually exceedingly difficult.

Let us denote by $B$ and $D$ the $3$-balls in which $X$ and $T$ lie, 
and let $S$ be the $4$-punctured sphere $S^2_4=\p B-X=\p D-T.$ For any $\phi$ in the mapping class group, $MCG(S^2_4),$ the action of $\phi$ on $S^2_4$ extends to maps of $B-X$ and of $D-T$ to $B-X'$ and $D-T'$, respectively, such that $N(X+T)=N(X'+T')$. In other words, $MCG(S^2_4)$ acts on all pairs $(X,T)$ while preserving $N(X+T)$. (The orbits of that action on tangles are the isotopy classes of tangles moving their endpoints around.)
We call pairs $(X,T)$ and $\phi(X,T),$ for $\phi\in MCG(S^2_4)$ \underline{equivalent}. 
Since all rational tangles belong to a single orbit of $MCG(S^2_4)$, this equivalence defines a $1$-$1$ correspondences between solutions of tangle equations \eqref{e-tangle} for different $T$'s.

Recall that a tangle is $2$-string if it has no loops.

\bpro\label{p.solnum}
(1) 
If $L$ is the unlink of 2 components then Eq. \eqref{e-tangle} has infinitely many solutions but only one $2$-string one.\\
(2) If $L$ is a knot or a non-trivial link then Eq. \eqref{e-tangle} has infinitely many $2$-string solutions.
\epro

\bpr
Since every tangle equation is equivalent to one with $T=\la 0\ra,$ we assume that for the proof. Hence, clearly it has at least one solution, $X$. 
For $L$ the unlink, $X=\la 0\ra$.

Assume that $L$ is not the unlink now. For any long knot $K$, let $Sat_X(K)$ be its 
satellite with $X$ as a companion tangle, as in Fig. \ref{f.satellite}.
\begin{figure}[h]
\centerline{\diag{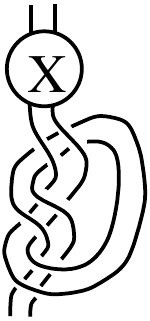}{0.4in}}
\caption{A satellite of a long trefoil}
\lb{f.satellite}
\end{figure}
Then $N(Sat_X(K))=N(X)$. Hence, if $X$ is a solution of Eq. \eqref{e-tangle}, then $Sat_X(K)$ is a solution as well. Since $Sat_X(K)$ and $Sat_X(K')$ are non-isotopic for different $K$ and $K'$ (here we assume that $X\ne \la 0\ra$), there are infinitely many different tangles of this form.
\epr 

For a nontrivial rational $L$ \cite[Thm. 2.2]{ES} shows infinitely many rational solutions to Eq. \eqref{e-tangle}.

Finally, we remark that Proposition \ref{p.solnum} does not hold it $T$ is allowed to be non-rational. For example, Eq. \eqref{e-tangle} has no solutions for certain pretzel tangles $T$ and $L$ the unknot, since one can show that some pretzel tangles do not embed into the unknot. 
(That fact is an easy application of Proposition \ref{p.embedKB}). 

\subsection{Systems of Tangle Equations}
\lb{s-systangleEq}

A system of tangle equations has the form \eqref{e-tangle0} or, equivalently,
$$\begin{cases}
 N(X+T_1)=L_1\\
 N(X+T_2)=L_2.
\end{cases}$$ 
where $T_1\ne T_2$.

Systems of tangle equations play a crucial role in the analysis of entanglement of DNA molecules.
The reason for that is that certain enzymes (called recombinase) separate circular DNA substrate molecules into two tangles: $T_1$, consisting of the part of the DNA molecule bound to the enzyme, and the other part, $X$, not bound to the enzyme. Then the enzyme replaces $T_1$ by a tangle $T_2$.
The substrate knot $L_1$ is controlled by the experiment. The link $L_2$ 
is observable in the experiment. Since tangles $T_1$ and $T_2$ are known,  DNA recombination processes lead to systems \eqref{e-tangle0}, where $X$ is the unknown.
See e.g. \cite{ADV, BM, Da2, DS, ES, Su, SECS, VLNSLDL,Ya} for further discussion of such systems.

The action $MCG(S^2_4)$ on tangles of the previous section extends onto an action on triples $(X,T_1,T_2)$. 
By analogy, we call triples $(X,T_1,T_2)$ and $\phi(X,T_1,T_2)$ \underline{equivalent} for any $\phi\in MCG(S^2_4).$
This equivalence defines a $1$-$1$ correspondences between solutions of systems of tangle equations \eqref{e-tangle0} for equivalent pairs $(T_1,T_2)$ and $(T_1',T_2')$.

\bpro\label{p-n-sols}
For every pair $(T_1,T_2)$ not equivalent to $(\la 0 \ra, \la k\ra)$ for $k\in \Z\cup \{\infty\}$ and for every $n>0$ there are  links $L_1$ and $L_2$ such that \eqref{e-tangle0} has at least $n$ solutions. 
\epro

A \underline{Montesinos tangle} is a sum of rational ones, $R_1+ \ldots +R_k$.
The numerator closure of a Montesinos tangle is a \underline{Montesinos link}.
In the above proposition, one can further assume that $L_1, L_2$ and the solutions $X$ are Montesinos.

\bpr[Proof of Proposition \ref{p-n-sols}]
By utilizing an equivalence, we can assume that Equations \eqref{e-tangle0} have form:
$$N(X)=L_1,\quad N(X+T)=L_2,$$
for some $L_1,L_2$ and $T=\la x_0\ra,$  for some $x_0\in \wh \Q -\Z.$ 
Let
$$X_i=T+...+T+\la x_1\ra+T+...+T,$$ 
for some $x_1\in \wh \Q,$ 
where there are $i$ $T$'s on the left and $2n-i$ $T$'s on the right.
By a cyclic symmetry, $N(X_i)$ coincide for $i=0, ..., n-1$ and also $N(X_i+T)$ coincide for $i=0, ..., n-1$.
 
We claim that choosing $x_1\in \wh\Q$ so that $x_1\ne x_0$ mod $1$ guarantees that all $X_1,...,X_n$ are distinct. Indeed, the links $N(X_i+\la x_0\ra)$ are Montesinos.  By the classification of Montesinos links, \cite{Bo}, they are all distinct if $x_0\ne x_1$ mod $1.$
\epr

In the case of systems of tangle equations, the conditions for the existence and uniqueness or finiteness of solutions 
seem unknown in general. 
We propose

\bcon
Every system (\ref{e-tangle0}) have finitely many solutions $X$ only.
\econ

We are going to see later (eg. in Sec. \ref{s.uniqueKB}) that the uniqueness of solutions is easier to analyze in the context of framed tangles and links. 

\section{Surgery Methods, Cosmetic Surgery Conjecture}
\label{s.surgery}

We denote by $\Sigma(T)$ the double cover of the $3$-ball $B^3$ branched along a tangle $T$ in it. 
Similarly, $\Sigma(L)$ denotes the double cover of $S^3$ branched along a link $L$. 
Since the double branched cover of a rational tangle $T$ is a solid torus, every solution $X$ of \eqref{e-tangle0} defines a knot (the core of $\Sigma(T)$) in $\Sigma(L).$

As an application of surgery methods we obtain 
Proposition \ref{p.LL-i} which we recall here:

\begin{named}{Proposition \ref{p.LL-i}}
\pLL
\end{named}

\bpr
Suppose that $X$ is an algebraic tangle solution to 
\beq\label{e.LL}
N(X+T_1)=L=N(X+T_2),
\eeq
for $T_1\ne T_2.$
By the above discussion, there is a non-trivial surgery of $\Sigma(L)$ yielding $\Sigma(L)$. If $L$ is rational then $\Sigma(L)$ is a lens space and if $L$ is a torus knot then $\Sigma(L)$ is Seifert with base the 2-sphere and three exceptional fibers, see eg. \cite[Sec. 3.1]{JP}. Either way it is irreducible and atoroidal. Furthermore, $\Sigma(X)$ is a graph manifold and, hence, non-hyperbolic. 
Hence, \cite[Cor. 1.4]{Ma} implies that $\Sigma(T_i)$ lies in a $3$-ball. Since $\Sigma(L)\ne S^3,$ the boundary of that ball is incompressible in $\Sigma(X),$ making it reducible. Since $X$ is rational or prime, that is impossible by \cite[Thm. 5]{Li}.
\epr

As observed in the proof above, every solution of \eqref{e.LL} yields a cosmetic surgery on $\Sigma(L)$. Such surgeries are rare and are subject of the Cosmetic Surgery Conjecture, \cite{Go, Ki}.  A special case of that conjecture 
is the Nugatory Crossing Conjecture, asserting that if a change of a sign of a crossing $c$ in a knot diagram $D$ results in an isotopic knot then $c$ is nugatory, that is $K$ has form below,
\cite{BFKP, BK, Ka, LM, To}. 
\begin{figure}[h]
\centerline{\diag{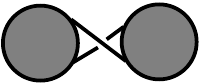}{.6in}}
\caption{A nugatory crossing in knot. Disks denote $1$-tangles.}
\lb{f-nugatory}
\end{figure}

We say that a tangle $X\subset B^3$ is \underline{unlinked} if it is a $2$-string tangle whose components can be isotoped so that they do not intersect each other.
(They can be all obtained by knotting the arcs of $\la 0\ra$ or of $\la \infty\ra$.) The following is straightforward:

\bpro
Nugatory Crossing Conjecture is equivalent to the following statement:
If $N(X+\la 1\ra)=N(X+\la -1\ra)$ is a knot then $X$ is unlinked.
\epro

\bpr 
(1) $\Rightarrow$ Suppose the Nugatory Crossing Conjecture holds and 
$$N(X+\la 1\ra)=N(X+ \la -1\ra)$$ 
is a knot $K$.
Then the crossing in $\la \pm 1\ra$ is nugatory in $K$. By the Nugatory Crossing Conjecture, $K$ consists of two arcs, each in a disjoint ball, as in Fig. \ref{f-nugatory}. Since these balls can be placed in the $3$-disk $B^3$ to form $X$, the statement follows.

$\Leftarrow$ Any knot with a nugatory crossing can be realized as $N(X+\la 1\ra)$ where the nugatory crossing is in $\la 1\ra.$ Since $X$ is unlinked the statement follows.
\epr

It is worth noting that $N(X+\la 1\ra)=N(X+\la -1\ra)$ cannot be a link for a $2$-string $X$.

Finally, the following example shows that the assumption of non-triviality of $L_1=L_2$ in Proposition \ref{p.LL-i} is essential. It will be useful later.

\blem\lb{l-unframed}
$X=\la \infty\ra$ and $\la -1/2\ra$ are the only 
solutions to
\beq\lb{e-taneq-eg}
N(X)=U=N(X+\la 1\ra),
\eeq
where $U$ denotes the unknot.
\elem

\bpr 
The numerator closure $N(\cdot)$ and the operation $N(\cdot +\la 1\ra)$, on the level of the double covers correspond to two different Dehn fillings of $\Sigma(X).$ Each of them yields the double cover of $S^3$ branched along $U$, i.e. $S^3.$ By a theorem of Gordon and Luecke, \cite[Thm. 2]{GL}, $\Sigma(X)$ must be a solid torus and, hence, $X$ is rational, $X=\la p/q\ra.$

By a theorem of Schubert, 
$$N(\la p/q\ra)=U=N(\la 1\ra)$$ 
only if $p=1$, \cite{KL, Sch}.
Since $N(\la\frac{1}{q}\ra+\la 1\ra)=U$, it is easy to see that $q$ is either $0$ or $-1/2.$
\epr

\section{Framed Tangles and Systems of Equations}
\label{ss.framedtangles}

\underline{Framed links} are tame embeddings of annuli $S^1\times I\cup ...\cup S^1\times I$ into $\R^3$, where $I$ is an open interval. Similarly, \underline{framed} \underline{tangles} are proper tame embeddings 
$$J_1\times I \cup J_2\times I \hookrightarrow D^2\times (-1,1),$$
where $J_1, J_2$ are closed intervals and the end arcs 
$\p J_1\times I$, $\p J_2\times I$ lie disjointly in $\p D^2=\p (D^2\times \{0\})$ each containing a different point from among $NE, SE, SW, NW$ in $\p D^2$.
Clearly, every link diagram and tangle diagram defines a framed link or tangle with its framing parallel to the page.  We require that every framed link and tangle can be represented in that way. (Hence, components with a half-twist framing are not allowed).

For our purposes it will be convenient to consider framed links and tangles up to \underline{balanced isotopy} given by the balanced Reidemeister moves, presented in Fig. \ref{f-balancedisotopy}.
\begin{figure}[h]
\centerline{\diag{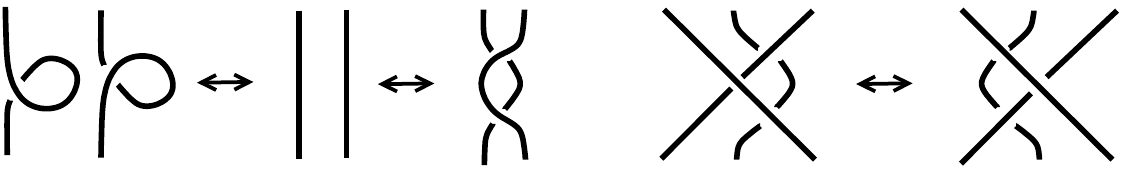}{3in}}
\caption{Balanced Reidemeister moves. (Diagrams have blackboard framing.)}
\lb{f-balancedisotopy}
\end{figure}
Note that this is a somewhat more flexible isotopy than the standard isotopy of framed links or tangles. 

A \underline{framed rational tangle} is a rational tangle with an arbitrary framing.
Given a diagram $D$ of a tangle or a link, we will denote by $D^n$ the framed diagram obtained from the page framing of $D$ by adding $|n|$ positive or negative twists, depending on the sign of $n$. Note that location of these twists does not matter up to balanced Reidemeister moves.

Definitions of tangle addition and of numerator and denumerator closures generalize immediately to framed tangles. 
Consequently, systems (\ref{e-tangle0}) can be considered in the context of framed tangles and links.

\subsection{Systems of Framed Tangle Equations.}

Framed links and tangles seem more appropriate for modeling DNA molecules which are double stranded.
Additionally, we will see soon, that it sometimes easier to solve tangle equations in the framed setting and framed solutions inform on unframed ones.

Let us assume that $T_1, T_2$  are framed rational tangles now. As before, we assume that $T_1, T_2$ are unequal as unframed tangles. (Note that if $T_1$ and $T_2$ in (\ref{e-tangle0}) differ by framing twists only, then by adjusting the framing of $L_2$, one can reduce (\ref{e-tangle0}) to a form in which $T_1=T_2$ as framed tangles. In that form, (\ref{e-tangle0}) is either inconsistent or reduces to a single equation.)

Note that Proposition \ref{p.solnum} easily extends to framed equations.

Given an unframed system (\ref{e-tangle0}) its \underline{framing} is a choice of a framing for $L_1,L_2,T_1,T_2.$
Let us consider the following framed version of system \eqref{e-taneq-eg} which will be useful later:
\beq\lb{e-tangle3}
N(X)=U^n\ \text{and}\ N(X+\la 1\ra)= U^m,
\eeq
for some $n,m\in \Z,$ where as defined above, $U^n$ denotes the unknot with $|n|$ twists of framing, positive or negative, depending on the sign of $n$.

\blem\lb{l-framed}
 (1) If $m=n+1$ then 
 $\la \infty\ra^n$ is the only solution to (\ref{e-tangle3}).\\
 (2) If $m=n-3$ then $\la -1/2\ra^{n+2}$  is the only solution to (\ref{e-tangle3})\\
 (3) If $m-n$ is neither $1$ or $-3$ then (\ref{e-tangle3}) has no framed solution.
\elem

\bpr Suppose that $X$ satisfies the above equations. Then stripped of its framing, it is either 
$X=\la \infty\ra$ or $\la -1/2\ra$, by Lemma \ref{l-unframed}.
Since\\ 
\centerline{\diag{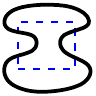}{0.4in} $=U^{f}$ and \diag{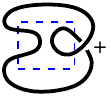}{0.4in} $=U^{f+1}$}
 where the dashed square marks the tangle $\la \infty\ra^f,$ for some $f\in \Z$, the tangle
 $X=\la \infty\ra$ lifts to a framed solution only iff $m-n=1$ and that solution is $\la \infty\ra^f$, where $f=n$.
 (For reader's convenience we marked the signs of crossings in the diagrams.)
Similarly,\\
\centerline{\diag{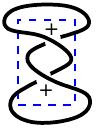}{0.4in} $=U^{f+2}$ and \diag{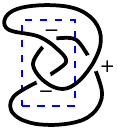}{0.5in} $=U^{f-1}$}
shows that $X=\la -1/2\ra$ lifts to a framed solution only iff $m-n=-3$ and that solution is $\la -1/2\ra^f$ for $f=n-2$.

For $m-n\ne 1, -3$, these equations are contradictory. 
\epr

Note that every framed system \eqref{e-tangle0} defines an unframed one and every framed solution descends to an unframed one. 
Lemmas \ref{l-unframed} and \ref{l-framed} show that not every solution to the unframed system can be lifted to a framed one. In fact, there may be no framing of \eqref{e-tangle0} for which all unframed solutions can be framed.

\section{The Kauffman bracket of a tangle}
\label{s-KB}

The Kauffman bracket $[L]\in \Z[A^{\pm 1}]$ is an invariant of framed unoriented links $L$ up to isotopy, satisfying the skein relations
$$\KPX = A \KPA +A^{-1}\KPB,\quad \KPC =\delta$$
(with the blackboard framing), where $\delta=-A^2-A^{-2},$
normalized so that the bracket of the trivially framed unknot $U^0$ is $[U^0]=1.$ 

By the above skein relations, each framed tangle can be expressed as 
$$T=[T]_0\cdot \KPB + [T]_\infty\cdot \KPA,$$ 
where $[T]_0,[T]_\infty\in \Z[A^{\pmo}]$ are uniquely defined. 
Recall that we call the vector
$[T]=\begin{pmatrix} [T]_0\\ [T]_\infty\end{pmatrix}$ the Kauffman bracket of $T.$ 

Kauffman brackets of tangles provide a criterium for embedding them into links:

The \underline{KB determinant} of $T$ denoted by $det_{KB}(T)\in \Q[A^\pmo]$ is $gcd([T]_0,[T]_\infty)$ in $\Q[A^\pmo]$. Note that it is well defined for unframed tangles up to a multiplicative factor of $A^n$, for $n\in \Z$.

\bpro\label{p.embedKB}
A necessary condition for a tangle $T$ embedding into a link $L$ is that $det_{KB}(T)$ divides $[L]$ in $\Q[A^\pmo]$.
\epro

This result generalizes the main result of \cite{Kr} which treats $A=e^{\pi i/4}$ case; see also \cite{Ru}.
(This idea was utilized for the trivial $L$ in \cite{TS1}.) More generally, a necessary condition is that $[L]\in \Z[A^\pmo]$ lies in the ideal $\la [T]_0,[T]_\infty\ra\triangleleft \Z[A^\pmo]$.\vs

\noindent {\it Proof of Proposition \ref{p.embedKB}:}
$$\TaddNumer{$T$}{$S$}\hspace*{-.18in}=[T]_0[S]_0 \diag{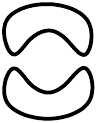}{0.32in} + [T]_0[S]_\infty\diag{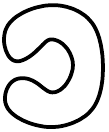}{0.32in} +[T]_\infty[S]_0\rotatebox[origin=c]{180}{\diag{0infty.pdf}{0.32in}}
+[T]_\infty[S]_\infty\KP{\filldraw[color=black, fill=none, thick] circle (0.23);
\filldraw[color=black, fill=none, thick] circle (0.45);}
$$
and, therefore, a tangle $T$ embeds into a link $L$ only if $gcd([T]_0,[T]_\infty)$ divides $[L]$ in $\Q[A^{\pm 1}].$ (Formally, one needs framed tangles and links for that, but since framing affects the bracket by a multiplicative factor of $(-A^3)^n$ which is a unit in $\Q[A^{\pm 1}]$, that statement makes sense in the unframed setting as well.)
\qed\vspace*{.1in}

Note that $[T]_0, [T]_\infty$ are integers for that $A=e^{\pi i/4}$.
The topological meaning of $gcd( [T]_0, [T]_\infty)$ for that $A$ seems to be an interesting question.
It appears to be related to the order of the torsion of the double branched cover of $T$,
however it is not always equal to it, as observed in \cite{Ru}.

The following discussion will be useful later.
We say that an unframed $2$-string tangle $T$ is of \underline{type $0$, $\pm 1$, or $\infty$} if it can be reduced to one of these rational tangles through crossing changes. Note that each tangle is of precisely one of these three types, depending on connections between its endpoints. 
Type $\pm 1$ has two subtypes: $+1$ and $-1$, depending on whether one can get from the original tangle $T$ to tangle $+1$ in an even or odd number of crossing changes between different components. 
For the statement below we will need to consider oriented tangles. The orientation will be considered up to the total orientation reversal, hence 
$$\KP{
   \filldraw[color=black, fill=none, thick] (0,0) circle (0.41);
    \draw[color=black, thick, <-] (-0.3,-0.3) -- (0.3,0.3);
    \draw[color=black, thick] (-0.3,0.3) -- (-0.05,0.05);
    \draw[color=black, thick, ->] (0.05,-0.05) -- (0.3,-0.3);
  }=\KP{
   \filldraw[color=black, fill=none, thick] (0,0) circle (0.41);
    \draw[color=black, thick, ->] (-0.3,-0.3) -- (0.3,0.3);
    \draw[color=black, thick, <-] (-0.3,0.3) -- (-0.05,0.05);
    \draw[color=black, thick] (0.05,-0.05) -- (0.3,-0.3);
  }.$$
Every oriented tangle of type $\pm 1$ is either oriented \underline{vertically} or \underline{horizontally} depending on whether one can reduce it to the vertically or horizontally
oriented tangles of Fig. \ref{f-verhor}, by crossing changes which involve an even number of changes between different components. 

If we denote the writhe of a framed oriented link $L$ by $w(L),$ then $(-A)^{-3w(L)}[L]$ is the Jones polynomial $J(L)$ for $t=A^{-4}.$ Since $J(L)\in t^{(|L|-1)/2}\Z[t^\pmo]$, where $|L|$ denotes the number of components of $L$, the bracket $[L]$ can be written with all exponents of $A$ congruent to 
$$3w(L) +2(|L|-1) \mod 4.$$
We extend this property to Kauffman brackets of tangles. 

\begin{figure}[h]
\KP{
   \filldraw[color=black, fill=none, thick] (0,0) circle (0.41);
    \draw[color=black, thick, <-] (-0.3,-0.3) -- (0.3,0.3);
    \draw[color=black, thick] (-0.3,0.3) -- (-0.05,0.05);
    \draw[color=black, thick, ->] (0.05,-0.05) -- (0.3,-0.3);
  }\quad
\KP{
   \filldraw[color=black, fill=none, thick] (0,0) circle (0.41);
    \draw[color=black, thick, ->] (-0.3,0.3) -- (0.3,-0.3);
    \draw[color=black, thick, <-] (-0.3,-0.3) -- (-0.05,-0.05);
    \draw[color=black, thick] (0.05,0.05) -- (0.3,0.3);
   }\hspace*{.3in}
 \KP{
   \filldraw[color=black, fill=none, thick] (0,0) circle (0.41);
    \draw[color=black, thick, ->] (-0.3,-0.3) -- (0.3,0.3);
    \draw[color=black, thick] (-0.3,0.3) -- (-0.05,0.05);
    \draw[color=black, thick, ->] (0.05,-0.05) -- (0.3,-0.3);
  }\quad
\KP{
   \filldraw[color=black, fill=none, thick] (0,0) circle (0.41);
    \draw[color=black, thick, ->] (-0.3,0.3) -- (0.3,-0.3);
    \draw[color=black, thick] (-0.3,-0.3) -- (-0.05,-0.05);
    \draw[color=black, thick, ->] (0.05,0.05) -- (0.3,0.3);   
     }
\caption{$1$ and $-1$ tangle types oriented vertically (on the left) and horizontally (on the right).} 
\lb{f-verhor}
\end{figure}

Let $|T|$ denote the number of connected components of a tangle $T.$

\bpro\label{p.exp}
For any oriented tangle $T$, the bracket $[T]_\tau$, for $\tau$ is $0$ or $\infty$, can be written with all exponents of $A$ congruent to 
$$3w(T)+2 |T|+2\eta(T) +\delta_\tau\mod 4$$ 
where $\delta_\tau =0$ or $2$ depending on whether $\tau$ is $0$ or $\infty$ and
$\eta(T)=0$ for $T$ of type $\pm 1$ vertically oriented (see Fig. \ref{f-verhor}) and for type $0$. Otherwise, $\eta(T)=1$; that is for type $\pm 1$ oriented horizontally and for $T$ of type $\infty$.
\epro

\bpr We prove it for all oriented tangle diagrams $D$ by induction on the crossing number $c(D)$ of $D$. The statement is obvious for $c(D)=0.$ Assume that the statement holds for all diagrams with $c(D)\leq C$. We are going to prove the statement for diagrams $D$ with $C+1$ crossings by the induction on the \underline{undescendedness} $ud(D)$, of $D$, which measures the failure of being a descending diagram.

It is defined as follows: We will refer to the connected components of the tangle represented by $D$ as the components of $D$ for simplicity. Let us also orient and order all these components and let us choose a base point on each loop component of $D$ away from crossings.  Next let us follow all components of $D$ according to their order, and according to their orientations, starting with their beginning (if they are arcs) or their base points (if they are loops). A crossing of $D$ is \underline{improper} if its underpass is transversed before its overpass is. (The improperness of a crossing will depend on the above choices.)
The undescendedness, $ud(D)$, of $D$ is the minimum of the numbers of improper crossings in$D$ over all choices of basepoints and choices of orientations and orderings of the components of the tangle represented by $D$. In particular, 
$$0 \leq ud(D)\leq c(D).$$

If $ud(D)=0$ then $D$ is isotopic to $\la 0\ra, \la \pm 1\ra$ or $\la \infty\ra$ depending on the type of the tangle represented by $D$, with possibly some additional unlinked trivial loop components. It is easy to see that the statement follows in this case. Assume that the statement holds for all diagrams $D$ with $c(D)=C+1$ and $ud(D)\leq u$. Let $D$ be a diagram with $c(D)=C+1$ and $ud(D)= u+1$ and let $v$ be one of the improper crossings counted in $ud(D)$ (according to a certain choice of base points, orientations, and of ordering). Assume that this crossing has sign $\ve=\pm$. Then 
$$A [D_+] - A^{-1} [D_-] =(A^2-A^{-2})[D_0],$$
where $D_\ve$ is $D$, and $D_{-\ve}$ is $D$ with the crossing $v$ reversed and $D_0$ is $D$ with $v$ smoothed out. Consequently,
$$[D] = A^{-2\ve}[D_{-\ve}] +\ve\cdot A^{-\ve}(A^2-A^{-2})[D_0].$$  
Since 
$$ud(D_{-\ve})<ud(D), \quad c(D_0)< c(D),$$
and $$w(D_{-\ve})=w(D)-2\ve \ \text{and}\ w(D_0)=w(D)-\ve,$$
by the inductive assumption,
the exponents of $A^{-2\ve}[D_{-\ve}]_\tau$ are
$$-2\ve+3(w(D)-2\ve) +2 |D| +2\eta(D)+\delta_\tau \mod 4$$
which equals 
\beq\label{e.exp1}
3 w(D) +2|D| +2\eta(D)+\delta_\tau \mod 4.
\eeq 
Similarly, the exponents of $A^{-\ve}(A^2-A^{-2})[D_0]_\tau$ are 
$$-\ve\pm 2 +3(w(D)-\ve)+2 |D_0| +2\eta(D_0)+\delta_\tau \mod 4$$
which equals
\beq\label{e.exp2} 
\pm 2 +3 w(D) +2|D_0|+2\eta(D_0) +\delta_\tau\mod 4.
\eeq

To complete the argument, let us assume first that $v$ is a single component crossing.
Then $D,D_{-\ve}, D_0$ have the same type and 
$$|D_\ve|=|D_{-\ve}|=|D_0|\pm 1,$$
Hence, the quantities in  \eqref{e.exp1} and \eqref{e.exp2} coincide and imply that the exponents of $[D]_\tau$ are
$$3 w(D)+2 |D| +2\eta(D)+\delta_\tau \mod 4,$$
as claimed in the statement, thus completing the proof of the inductive step in the single-component crossing case.

If $v$ is a two-component crossing, then 
$D,D_{-\ve}$ have the same type and 
$$|D_\ve|=|D_{-\ve}|=|D_0|.$$
A two-component crossing smoothing transforms tangle types and orientations in one of the following three ways:\vs

\centerline{
\KP{
   \filldraw[color=black, fill=none, thick] (0,0) circle (0.41);
    \draw[color=black, thick, ->] (-0.3,0.3) -- (0.3,-0.3);
    \draw[color=black, thick, <-] (-0.3,-0.3)  -- (0.3,0.3);
    {}
  }\ $\leftrightarrow$\ 
 \KP{
  \filldraw[color=black, fill=none, thick] circle (0.41);
  \draw[color=black, thick, <-] (-0.3,-0.3) .. controls (0.02,0) .. (-0.3,0.3);
   \draw[color=black, thick, <-] (0.3,-0.3) .. controls (-0.02,0) .. (0.3,0.3);
    {}
  }, \quad 
\KP{
   \filldraw[color=black, fill=none, thick] (0,0) circle (0.41);
    \draw[color=black, thick, ->] (-0.3,0.3) -- (0.3,-0.3);
    \draw[color=black, thick, ->] (-0.3,-0.3)  -- (0.3,0.3);
  }\ $\leftrightarrow$\ 
\KP{
   \filldraw[color=black, fill=none, thick] (0,0) circle (0.41);
    \draw[color=black, thick, ->] (-0.3,0.3) .. controls (0,-0.02) .. (0.3,0.3);
    \draw[color=black, thick, ->] (-0.3,-0.3) .. controls (0,0.02) .. (0.3,-0.3);
  },\quad 
 \KP{
  \filldraw[color=black, fill=none, thick] circle (0.41);
  \draw[color=black, thick, <-] (-0.3,-0.3) .. controls (0.02,0) .. (-0.3,0.3);
   \draw[color=black, thick, ->] (0.3,-0.3) .. controls (-0.02,0) .. (0.3,0.3);
  } = 
  \KP{
  \filldraw[color=black, fill=none, thick] circle (0.41);
  \draw[color=black, thick, <-] (-0.3,-0.3) .. controls (0.3,0) .. (0.08,0.12);
   \draw[color=black, thick, ->] (0.3,-0.3) .. controls (-0.3,0) .. (0.3,0.3);
  \draw[color=black, thick] (-0.3,0.3) -- (-0.06,0.18);
  }\ $\leftrightarrow$\  
 \KP{
   \filldraw[color=black, fill=none, thick] (0,0) circle (0.41);
    \draw[color=black, thick, ->] (-0.3,0.3) .. controls (0,-0.02) .. (0.3,0.3);
    \draw[color=black, thick, <-] (-0.3,-0.3) .. controls (0,0.02) .. (0.3,-0.3);
  }.
  }\vs 
  (Since these are tangle types, no crossing signs are indicated.)\vs
  
  {}
Consequently, $\eta(D)\ne \eta(D_0).$
Hence, the quantities in  \eqref{e.exp1} and \eqref{e.exp2} coincide and imply that
the the exponents of $[D]_\tau$ are
$$3 w(D)+2 |D| +2\eta(D)+\delta_\tau \mod 4,$$
as claimed in the statement, thus completing the proof of the inductive step in the two-component crossing case.
\epr

\section{The KB ratio of a tangle}
\label{s-KBratio}
Let
$$Q(T)=[T]_\infty/[T]_0\in \Q(A)\cup \{\infty\}.$$ 
 It is well defined, because $[T]_\infty$ and $[T]_0$ cannot both vanish. (This follows for example from Proposition \ref{p.embedKB} and the fact that $[N(T)]\ne 0,$ because $[N(T)]=(-2)^{|N(T)|}$ for $A=\pm 1$.) 

Note that $Q(T)$ is preserved by the first Reidemeister move and, hence, it is an invariant of unframed tangles.

Let 
$$q= -A^4\ \text{and}\ \{ T\}_q=A^{-2}Q(T).$$

Proposition \ref{p.exp} implies
$$\{T\}_q\in \Q(q)\cup \{\infty\},$$
thus proving Theorem \ref{t.q-qi}.
We call it the \underline{Kauffman bracket ratio}, or the \underline{KB-ratio} of $T$, for short. 
KB-ratios will play an important role in the remainder of the paper.
 
\bpro\label{p.qrat}
For any tangle $T,$\\
(1) $\{-T\}_q=-q^{-1}\{ T\}_{q^{-1}},$ where the subscript on the right side indicates $q$ substituted by $q^{-1}$. (Recall that $-T$ is the mirror image of $T$.)\\
(2) $\{ R(T)\}_q= -q^{-1}/\{T \}_q$, where $R$ denotes the rotation operation, as before.

Additionally, for any $T'$,\\
(3) $\{T+T' \}_q=\{T\}_q+\{T' \}_q+(q-1)\{T\}_q\{T' \}_q.$
\epro

\bpr 
(1) $\{-T\}_q=A^{-2}Q(-T)=A^{-2}Q(T)_{A^{-1}}=A^{-4}\{T\}_{q^{-1}}.$\\
(2) Since $Q(R(T))=1/Q(T)$,
$$\{R(T)\}_q=A^{-2}Q(R(T))=A^{-2}/Q(T)=A^{-4}/(A^{-2}Q(T))=-q^{-1}/\{T\}_q.$$
(3) By the following identity
$$\Taddition{$T$}{$T'$} = [T]_0[T']_0 \KPB+ ([T]_0[T']_\infty+[T]_\infty[T']_0+ [T]_\infty[T']_\infty\delta) \KPA,$$
$$Q(T+T')=Q(T)+Q(T')+\delta Q(T)Q(T'),$$
and, hence,
$$\{T+T' \}_q=\{T\}_q+\{T' \}_q+A^2\delta\{T\}_q\{T' \}_q.$$
\epr

\section{$q$-deformed rationals}
\label{s-q-rationals}

Quantum integers,
$$[n]_q=\begin{cases} \frac{1-q^n}{1-q} 
& \text{for } n\ne 0\cr
0 & \text{for } n=0,\cr
\end{cases}$$
for $n\in \Z$, 
appear already in the 1808 work of Gauss on binomial coefficients, \cite{Gau} and in \cite{Ja}.
They are at the foundation of quantum calculus, \cite{KC}, and are indispensable in quantum algebra and in quantum topology, see eg. \cite{Kas}.

Quantum integers were extended in \cite{MO} to ``$q$-deformed rationals'' $[x]_q$ for $x\in \Q.$
They satisfy the following identities
 \beq\label{e.qrat}
 [x+1]_q=q[x]_q+1,\quad [-x]_q=-q^{-1}[x]_{q^{-1}},\quad [1/x]_q= \frac{1}{[x]_{q^{-1}}},
 \eeq
 cf. \cite{LMG}.
 {}

They can be defined also through our bracket $\{\cdot\}_q$ and, in fact, in this way they have been discovered independently by us in the first version of this paper on arXiv.

Note that the first two identities of Eqs. \ref{e.qrat} are satisfied by $\{\cdot\}_q$ by Proposition \ref{p.qrat}. 
The last one is satisfied as well: According to Conway's correspondence $\la 1/x\ra=R(-\la x\ra)$ and, hence,
$$\{\la 1/x\ra\}_q=\{R(-\la x\ra)\}_q=-q^{-1}/\{-\la x\ra\}_q=-q^{-1}/(-q^{-1}\{\la x\ra\}_q)_{q\to q^{-1}}=1/\{\la x\ra\}_{q^{-1}}.$$ 
By the theory of continued fractions, the above identities determine the values of $[x]_q$ and of 
$\{x\}_q$ for all $x\in \Q.$ (In fact, the formulas for $x+1$ and $-1/x$ are sufficient.) This implies the following statement from the Introduction:

\begin{named}{Theorem \ref{t.qrati}}
For every $x\in \wh \Q,$
$[x]_q=\{\la x\ra\}_q$.
\end{named}

Let $\{\la x\ra \}_1$ denote $\{\la x\ra \}_q$ evaluated at $q=1.$
By the above theorem and Eq. \eqref{e.qrat}, we have:
\bcor[\cite{MO}]\lb{c-qi-1-1}
For every $x\in \wh\Q,$ $\{\la x\ra \}_1=x$. 
\ecor

\bcor\label{c.ratuniq}
The Kauffman bracket distinguishes all framed tangles among rational ones.
\ecor

\bpr Any two framed rational tangles $T,T'$ with $\{T\}_q=\{T'\}_q$ must represent the same rational number and, hence, may differ by framing only. However, if $T'$ is obtained from $T$ by adding $n$ twists to it (where negative $n$ means $|n|$ negative twists) then $[T']= (-A)^{3n}[T].$ Hence, $n=0.$
\epr

\section{The uniqueness of the KB of a solution of a system of tangle equations}
\label{s.uniqueKB}

\bthm\label{t.uniqueKB}
For any system of framed tangle equations (\ref{e-tangle0}), $[X]$ is unique.
\ethm

\bpr
Let $X$ be a solution of framed (\ref{e-tangle0}). 
Then 
$$\TaddNumer{$X$}{$T_1$}\hspace*{-.18in}=[X]_0[T_1]_0 \diag{00.pdf}{0.32in} + [X]_0[T_1]_\infty\diag{0infty.pdf}{0.32in} +[X]_\infty[T_1]_0\rotatebox[origin=c]{180}{\diag{0infty.pdf}{0.32in}}
+[X]_\infty[T_1]_\infty\KP{\filldraw[color=black, fill=none, thick] circle (0.23);
\filldraw[color=black, fill=none, thick] circle (0.45);},
$$
and analogously for $T_2.$ Hence,
$$\left(\bmat [N(X+T_1)]\\ [N(X+T_2)] \emat\right)=B\cdot \left(\bmat [X]_0\\ [X]_\infty \emat\right)$$
where
$$B=\left(\bmat [T_1]_0\delta + [T_1]_\infty & [T_1]_0 + [T_1]_\infty\delta\\ 
[T_2]_0\delta + [T_2]_\infty & [T_2]_0 + [T_2]_\infty\delta\\
\emat\right).$$
Note that
$$\det B=(\delta^2-1)\cdot \det \left(\bmat [T_1]_0 & [T_1]_\infty\\ [T_2]_0 & [T_2]_\infty\emat\right)$$
is non-zero, even for $q=1$ by Corollary \ref{c-qi-1-1}.

Given a system of equations (\ref{e-tangle0}),
let $$q_\mu=\frac{\det \left(\bmat [L_1] &  [T_1]_0+\mu [T_1]_\infty\\ [L_2] &  [T_2]_0+\mu [T_2]_\infty\emat\right)}
{(\delta-\mu)\cdot \det \left(\bmat [T_1]_0 & [T_1]_\infty\\ [T_2]_0 & [T_2]_\infty\emat\right)}\ \text{for}\ \mu=\pm 1.$$

Then, by Cramer's Rule, one can verify that
$$[X]_0=\det \left(\bmat [L_1]&  [T_1]_0+[T_1]_\infty \delta\\ [L_2]&  [T_2]_0+ [T_2]_\infty\delta\emat\right)/\det B=
\frac{1}{2}(q_1-q_{-1})$$
and 
$$[X]_\infty=\det \left(\bmat [T_1]_0\delta+[T_1]_\infty & [L_1]\\  [T_2]_0\delta+ [T_2]_\infty & [L_2]\emat\right)/\det B= \frac{1}{2}(q_{1}+q_{-1}).$$
\epr

The above equations provide necessary algebraic conditions for the existence of a framed solution of 
(\ref{e-tangle0}). In particular, we have

\bcor\lb{c-necessary}
A necessary condition for the existence of a framed solution to (\ref{e-tangle0}) is that  
$$p_0 :=\frac{1}{2}(q_1- q_{-1}),\quad  p_\infty :=\frac{1}{2}(q_1+ q_{-1})\ \in \Z[A^{\pm 1}]$$ 
and that $A^{-2}p_\infty/p_0\in \Q(A^4).$
\ecor

Note that Theorem \ref{t.uniqueKB} together with Corollary \ref{c.ratuniq} implies Proposition \ref{p.rat-uniqe-i} stating that every framed system (\ref{e-tangle0}) has at most one framed rational solution. 

\section{Jones Conjecture for Tangles}
\lb{s-JC-for-tan}

Recall that $U^n$ denotes the unknot with framing $n\in \Z$.

\blem[Kauffman bracket version of Jones Conjecture]\lb{l-KBJC}\ \\ 
The Jones conjecture (JC) is equivalent to its Kauffman bracket version (KBJC): 
 if $[K]=r\cdot A^k$, for some $r,k\in \Z$ then $K=U^n$  for some $n\in \Z.$ 
\elem

\bpr KBJC $\Rightarrow$ JC: Suppose that the Jones polynomial 
of $K$ is $J(K)=1$ for some knot $K$. Then $K$ with some framing
has its Kauffman bracket equal to $(-A^3)^n$ for some $n\in \Z.$ By KBJC, $K=U^n$. Hence, $K$ is trivial as an unframed knot.

JC $\Rightarrow$ KBJC: Suppose that $[K]=r\cdot A^k$, for some $r,k\in \Z$. Then $J(K)=r (-A)^{-3w(K)}A^k.$ By \cite[Cor. 3]{Gan},
$J(K)=1$. Hence, $K$ is (unframed) trivial, by the Jones conjecture. 
\epr

Now we can formulate three versions of the Jones conjecture for tangles:

\bthm\lb{JC=JCT}
The Jones conjecture is equivalent to each of the following statements:\\
(a) For any framed tangle $T$, if $[T]=\begin{pmatrix} 1\\ 0\end{pmatrix}$
then $T=\la 0\ra^0.$\\ 
(b) If $[T]=\begin{pmatrix} r\cdot A^{n}\\ 0\end{pmatrix}$
for some $r,n\in \Z$ then $T=\la 0\ra$ as an unframed tangle.\\
(c) If $[T']=[T]$ and $T$ is rational then $T'=T$ as framed tangles. 
\ethm

It is woth noting that statement (c) above does not hold for rational knots. That is there are examples of distinct rational knots with coinciding Jones polynomials, \cite{Kan}.\vspace*{.1in}

\noindent{\it Proof of Theorem \ref{JC=JCT}:}  JC implies (a): Assume the Jones conjecture holds and that $[T]=(1,0).$ (For convenience, in this proof we will write all vertical vectors horizontally.)
Then $[R(T)]=(0,1)$ 
and it is easy to check that
$$[N(R(T))]=1\ \text{and}\ [N(R(T)+\la 1\ra)]=-A^3.$$ 
By Lemma \ref{l-KBJC}, 
$$N(R(T))=U^0\ \text{and}\ N(R(T)+\la 1\ra)=U^1.$$ 
Now, by Lemma \ref{l-framed}, $R(T)=\la \infty\ra$ and, hence, $T=\la 0\ra.$

(a) implies (c): Suppose that $T$ is rational and $[T']=[T]$.
Then $T$ can be transformed 
$$T=T_1\to ... \to T_k=\la 0\ra$$ 
by the operations of rotation, $R(\cdot),$ of addition of one, $P(T)= T +\la 1\ra,$ its inverse, $P^{-1}(T).$ and adding a positive or negative kink, $F^{\pm 1}(T)$, (for framing changes).
Let us apply the same operations to $T':$
\beq\label{e.T''}
 T'=T_1'\to ... \to T_k'=T''.
 \eeq 
 
Since
$$[T+\la 1\ra]=\left[\TwithCross\right]=[T]_0\left[\KPX\right]+[T]_\infty [\KPAKink]=
[T]_0\left(\bmat A^{-1}\\A \emat\right)+[T]_\infty\left(\bmat 0\\-A^3 \emat\right),$$
the operations $R$ and $P$ 
induce linear transformations of the Kauffman brackets: 
$$[R(T)] = \left(\bmat 0 &1\\ 1 & 0\emat\right)[T],\quad [P(T)] = [T+\la 1\ra]=\left(\bmat A^{-1} & 0\\ A & -A^3\emat\right)\cdot [T].$$
 
Since $R,P$ and $F$ are invertible, 
 $[T_i]=[T_i']$ for every $i$ and, hence, $[T'']=(1,0).$
By (a), $T_k'=\la 0\ra$, implying that $T'=T.$

(c) implies (b): Assume that $[T]=(r\cdot A^n,0)$, for some $r,n\in \Z.$ Then $[D(T)]=r\cdot A^n$ and by \cite[Cor 3]{Gan} (as in the proof of Lemma \ref{l-KBJC}), $J(D(T))=1$. That implies that $[T]=((-A^3)^k,0)$ has the bracket of $\la 0\ra$ with some framing $k\in \Z.$
Hence, by (c), $T=\la 0\ra$, as unframed tangle. 

(b) implies (a): Suppose $[T]=(1,0).$ Then by (b), $T=\la 0\ra$ as an unframed tangle.
Since $[T]=(1,0)$, the framing of $T$ must be trivial.

(a) implies JC: Let $K$ be a knot with trivial Jones polynomial. Let us frame it so that $[K]=1.$ Let $K\# \la 0\ra$ be the connected sum of $K$ with the lower strand of $\la 0\ra.$ Then  $[K\# \la 0\ra]=(1,0)$
and $K\#\la 0\ra=\la 0\ra$ as unframed tangles, by (b) (which we proved is implied by (a)), implying that $K$ is trivial.
\qed\vspace*{.1in}

Theorem \ref{JC=JCT} can be further refined to consider its statements up to a certain numbers of crossings.
The \underline{crossing number}, $c(T)$ of a framed or unframed tangle $T$ is the minimal number of crossings in its unframed isotopy class. Then we can strengthen the implication  JC $\Rightarrow$ (c) above as follows:

\bpro
Assume that the Jones conjecture holds for knots up to $N$ crossings. If $[T']=[T],$ and $c(T)+c(T')<N$ and 
$T$ is rational then $T'=T$ (as framed tangles). 
\epro

\noindent{\it Proof} follows the above proof for $JC \Rightarrow (a) \Rightarrow (c).$  Note that the tangle $T''$ in Eq. \eqref{e.T''} has at most $c(T)+c(T')$ crossings. Now observe that the Jones conjecture up to $N$ crossings implies statement (a) up to $N-1$ crossings.
\qed

Note that Theorems \ref{t.uniqueKB} and \ref{JC=JCT} imply Theorem \ref{JC-UTESi}.

\section{Slopes of Tangles}
\lb{s.slope}
Recall that we call $\{T\}_1\in \wh\Q$  the algebraic slope of $T$ and that we postulate in Conjecture \ref{c.slopeeq} that if an $m$-essential surface in the complement of $T$ has a slope, that slope is $\{T\}_1.$
In this section we define the slope of a surface and we prove our conjecture for algebraic tangles. 

Let us denote by $\cal T$ the double cover of the boundary $S^2$ of the $3$-ball containing $T$, branched around $T\cap S^2$ and let $\theta: \cal T\to \cal T$ be the Deck transformation of that covering. 
The arc $\mu$ of the denumerator closure arc connecting SW and NW and its image $\theta(\mu)$ under $\theta$ form a loop which we call the meridian of $\cal T$, cf. Fig. \ref{f-torus-cover}. (By considering the numerator arc SE-NE, we get a parallel loop in $\cal T.$) Similarly, the arc SW-SE of the numerator closure together with its $\theta$-image forms a loop which we call the longitude. We orient the meridian and longitude as in Fig. \ref{f-torus-cover}.

Recall that a loop $\gamma$ in $S^2_4$ is dividing if it separates the punctures of $S^2_4$ into two groups of two. Its \underline{slope} $s(\gamma)$ is the slope of its lift to $\cal T$ in the longitude-meridian basis of $\cal T,$ so that the slopes of the longitude and the meridian in $\cal T$ are of the slope $0$ and $\infty$, respectively. 
It is well defined by \cite[Prop. 2.6]{FM}.

\begin{figure}[h]
\centerline{\diag{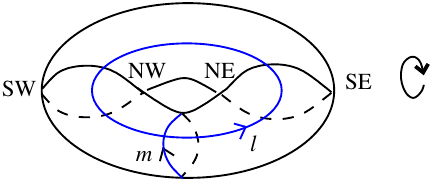}{2.3in}}
\caption{The double cover of a sphere $\p B^3$ branched along the endpoints of a tangle in $B^3$.}
\lb{f-torus-cover}
\end{figure}

As before, consider an $m$-essential surface $S$ properly embedded in $D^3 \smallsetminus T$. 
Recall that if its boundary contains a dividing loop $\gamma$ we say that $S$ is dividing and has slope $s(\gamma)$. (Then it is easy to see that all other dividing loops of $\p S$ must be parallel.) 
This definition is inspired by that in \cite{Oz2}, where the slopes are defined for algebraic tangles.


For example, $T=\la 0\ra$ has a horizontal disk $S=D^2$ in its complement which is $m$-essential. Since $\p D^2$ lifts to the longitude in $\cal T$, the slope of $\la 0\ra$ is $0$. One can construct any rational tangle from $\la 0\ra$ by the operations of rotation and addition of $\pm 1$ which modify $S$ accordingly. Consequently, it is easy to see that the slope of any rational tangle $\la x\ra$ is $x.$ 

Ozawa proves that for every algebraic tangle $T$, every $m$-essential surface in $B^3-T$ is dividing and its slope depends on $T$ only, \cite{Oz2}. This is the \underline{slope of $T$}, $s(T)\in \wh \Q.$ 

Recall Theorem \ref{t.slope-i}:

\begin{named}{Theorem \ref{t.slope-i}}
For every algebraic tangle $T$, $\{T\}_1=s(T).$
\end{named}

\bpr
By \cite[Proposition 2.4]{Oz2},
$$s(T_1+T_2)=s(T_1)+s(T_2),\quad s(-T)=-s(T),\quad s(R(T))= -1/s(T)$$
for  algebraic tangles. Furthermore, $s(\la n\ra)= n$, for $n\in \Z$. 
Since $\{\cdot\}_1$ satisfies the same properties by Proposition \ref{p.qrat} and 
since these properties determine the values of $s(T)$ and of $\{T\}_1$ for algebraic tangles, the statement follows.
\epr

The following provides an interpretation of the algebraic slope in terms of algebraic topology:

\bpro
$|\{T\}_1|= \det(N(T))/\det(D(T))$
where $\det$ denotes the link determinant. 
\epro

\bpr
The determinant of $N(T)$ is 
$$|J(N(T),-1)|=|(-A^3)^{-w(N(T))}[N(T)] |=| [N(T)] |,$$ 
where $A=e^{\pi i/4}$, $w(N(T))$ is the writhe of $N(T)$ (with some orientation) 
and 
$$[N(T)]=\diag{inftyN.pdf}{.35in}[T]_\infty+\diag{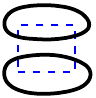}{.35in}[T]_0=
[T]_\infty.$$ 
Taking an analogous formula for $[D(T)]$ we have
$$\det(N(T))/\det(D(T))=|[T]_\infty/[T]_0|= |Q(T)_1|=|A^2\{T\}_1|=|\{T\}_1|.$$ 
\epr

%


\begin{thebibliography}{9999}
%



\bibitem[ADV]{ADV} J. Arsuaga, Y. Diao, M. Vazquez, Mathematical Methods in Dna Topology: Applications to Chromosome Organization and Site-Specific Recombination, in ``Mathematics of DNA Structure, Function and Interactions,''  The IMA Volumes in Mathematics and its Applications book series Vol. 150, Springer, 7--36.



\bibitem[BFKP]{BFKP} C. Balm, S. Friedl, E. Kalfagianni, and M. Powell, Cosmetic crossings and Seifert matrices, {\em Comm. Anal. Geom.} {\bf 20} (2012), no. 2, 235--253,

\bibitem[BK]{BK} C.J. Balm, E. Kalfagianni, Knots without cosmetic crossings, 
{\em Topology Appl.} {\bf 207} (2016), 33--42, arXiv:1406.1755


\bibitem[Bo]{Bo} F. Bonahon, Involution et fibres de Seifert dana lea vari\'et\'es de dimension 3," Th\'ese
de 3e cycle, Orsay, 1979.


\bibitem[BM]{BM} D. Buck, M. Mauricio, Connect sum of lens spaces surgeries: application to Hin recombination, {\em Math. Proc. Camb. Phil. Soc.}  {\bf 150} (2019), 505--525.

\bibitem[BZ]{BZ} G. Burde, H. Zieschang, Knots, de Gruyter Studies in Mathematics, 5 (1985).


\bibitem[Co]{Co} J.H. Conway, An enumeration of knots and links, in Computational Problems in Abstract Algebra, J. Leech (ed.), Pergamon Press (1969) 329--358. 



\bibitem[Da2]{Da2} I.K. Darcy, Modeling protein-DNA complexes with tangles, {\em Comput. Math. Appl.} {\bf 55}, no. 5, (2008),  924--937.

\bibitem[DS]{DS} I.K. Darcy, R.G. Scharein, TopoICE-R: 3D visualization modeling the topology of DNA recombination, 
{\em Bioinformatics,} {\bf 22} no. 14, (2006) 1790--1791.




\bibitem[EKT]{EKT} S. Eliahou, L.H. Kauffman, M.B. Thistlethwaite, Infinite families of links with trivial Jones polynomial, {\em Topology} {\bf 42} (2003) 155--169.



\bibitem[ES]{ES} C. Ernst, D.W. Sumners, A calculus for rational tangles, {\em Math. Proc. Camb. Phil. Soc.} {\bf 108} (1990) 489--515.  

\bibitem[FM]{FM} B. Farb, D. Margalit, Primer on Mapping Class Group, Princeton Mathematical Series, 41, 2011.


\bibitem[Gan]{Gan} S. Ganzell, Local moves and restrictions on the Jones polynomial, {\em J. Knot Theory and its Ramifications}, {\bf 23} (2014) 1450011.

\bibitem[Gau]{Gau} C.F. Gauss, Summatio quarumdam serierum singularium, 1808, Gottingae,
\url{http://resolver.sub.uni-goettingen.de/purl?PPN602151724}

\bibitem[GL]{GL} C.McA. Gordon, J. Luecke, Knots are determined by their complements, {\em Bull. Amer. Math. Soc.} {\bf 20}, 1, (1989) 83--87.

\bibitem[Go]{Go} C. Gordon, Dehn surgery on knots, Proceedings of the International
Congress of Mathematicians, Vol. I, II (Kyoto, 1990), 631–642, Math. Soc.
Japan, Tokyo, 1991. 


\bibitem[Ja]{Ja} F.H Jackson, On generalized functions of Legendre and Bessel,
{\em Trans. Roy. Soc. Edinburgh,} {\bf 41} (1903), 1--28

\bibitem[JP]{JP} Y. Jang, L. Paoluzzi, Double branched covers of tunnel number one
knots, {\em Geometriae Dedicata}  {\bf 211} (2019), 129--143, arXiv: 1905.05366.

\bibitem[Jo]{Jo} V.F.R. Jones, Ten problems, in Mathematics: Frontiers and Perspectives, American Mathematical Society Providence, (2000), 79--91.

\bibitem[KC]{KC} V. Kac and P. Cheung, Quantum calculus, Universitext, Springer-Verlag, New York,
2002.

\bibitem[Ka]{Ka} E. Kalfagianni, Cosmetic crossing changes of fibered knots, {\em J. Reine Angew.
Math.} {\bf 669} (2012), 151--164.

\bibitem[Kan]{Kan} T. Kanenobu, Examples of Polynomial Invariants of Knots and Links, II, {\em Osaka J. Math}, {\bf 26} (1989) 465--482

\bibitem[Kas]{Kas} C. Kassel, Quantum groups, Graduate Texts in Mathematics, vol. 155, Springer-Verlag,
New York, 1995.

\bibitem[KL]{KL} L.H. Kauffman, S. Lambropoulou, 
Classifying and Applying Rational Knots and Rational Tangles, in Physical Knots: Knotting, Linking and Folding Geometric Objects in $\R^3$ (Proc. Conference Las Vegas 2001), J.A. Calvo, et al, eds, Contemporary Math, {\bf 304} American Math. Soc. 2002, 223--259. 


\bibitem[Ki]{Ki} R. Kirby (ed.), Problems in low-dimensional topology, in: Geometric topology (Athens, GA 1993), AMS/IP Stud. Adv. Math. 2 (part 2), American Mathematical Society, Providence (1997), 35--473, \url{https://math.berkeley.edu/~kirby}

\bibitem[Kr]{Kr} D. Krebes, An obstruction to embedding $4$-tangles in links, {\em J. Knot Theory Ramif.} {\bf 8} (1999), 321--352, arXiv:9902119



\bibitem[LR]{LR} R. Lawrence, O. Rosenstein, Jones rational coincidences, arXiv: 2105.13897


\bibitem[LMG]{LMG} L. Leclere, S. Morier-Genoud, q-Deformations Of The Modular Group And Of The Real Quadratic Irrational Numbers, {\em Advances in Applied Mathematics}
{\bf 130} September 2021, 102223, arXiv: 2101.02953

\bibitem[Li]{Li} W.B.R. Lickorish, Prime knots and tangles, {\em Trans. Amer. Math. Society}, {\bf 267}, (1981) 321--332.

\bibitem[LM]{LM} T. Lidman, A.H. Moore, Cosmetic Surgery In L-Spaces and Nugatory Crossings,
{\em Trans. of A.M.S} {\bf 369} No. 5, (2017), 3639--3654.


\bibitem[MO]{MO} S. Morier-Genoud, V. Ovsienko, q-deformed rationals and q-continued fractions, 
{\em Forum of Mathematics, Sigma,} {\bf 8} (2020) E13  doi:10.1017/fms.2020.9
arXiv:1812.00170.


\bibitem[Ma]{Ma} D. Matignon, On the knot complement problem for non-hyperbolic knots, {\em Topology and its Applications}, {\bf 157} (2010) 1900--1925.


\bibitem[Mo2]{Mo2} J.M. Montesinos, Revetements ramifies des noeuds, Espaces fibres de Seifert et scindements de Heegaard, {\em Publicaciones del Seminario Mathematico Garcia de Galdeano} Serie II, Seccion 3 (1984).



\bibitem[Oz1]{Oz1} M. Ozawa, Morse position of knots and closed incompressible surfaces, {\em J. Knot Theory and its Ramifications}, {\bf 17} (2008) 377--397.

\bibitem[Oz2]{Oz2} M. Ozawa, Rational structure on algebraic tangles and closed incompressible surfaces in the complements of algebraically alternating knots and links, {\em Topology and its Applications} {\bf 157} no. 12, (2010) 1937--1948.


\bibitem[Ru]{Ru} D. Ruberman, Embedding tangles in links, 
{\em J. Knot Theory Ramif.} {\bf 9} No. 4 (2000), 523--530, arXiv: 0001141.

\bibitem[Sch]{Sch} H. Schubert, Knoten mit zwei Br\"ucken, {\em Math. Zeitschrift,} {\bf 65} (1956), 133--170.




\bibitem[Su]{Su} D.W. Sumners, Untangling DNA, {\em Math. Inteligencer}, {\bf 12} No. 3 (1994) 71--80.
\bibitem[SECS]{SECS} D.W. Sumners, C. Ernst, N.R. Cozzarelli, S.J. Spengler, Mathematical analysis of the mechanisms of DNA recombination using tangles, {\em Quart. Rev. of Biophys.}, {\bf 28} (1995) 253--313.


\bibitem[To]{To} I. Torisu, On nugatory crossings for knots, {\em Topology Appl.} {\bf 92} (1999), no. 2, 119--129.

\bibitem[TS1]{TS1} R. E. Tuzun, A. S. Sikora, Verification of the Jones unknot
conjecture to 22 crossings, {\em J. Knot Theory and its Ramifications}
{\bf 27} (2018) 1840009,  arXiv: 1606.06671.

\bibitem[TS2]{TS2} R. E. Tuzun, A. S. Sikora, Verification of the Jones unknot
conjecture to 24 crossings, {\em J. Knot Theory and Its Ramifications} {\bf 30} no. 3 (2021) {\bf 27}, 2150020, arXiv: 2003.06724. 

\bibitem[VLNSLDL]{VLNSLDL} A.A. Vetcher, A.Y. Lushnikov, J. Navarra-Madsen, R.G. Scharein, Y.L. Lyubchenko, I.K. Darcy, S.D. Levene, DNA Topology and Geometry in Flp and Cre Recombination, {\em Journal of Molecular Biology,} {\bf 357}(4) (2006) 1089--1104.




\bibitem[Ya]{Ya} J. Yang, Non-hyperbolic solutions to tangle equations involving composite links, arXiv: 1709.01785

\end{thebibliography}
\end{document}